# Resilience Enhancement of Electricity-Computing Interdependency System against Contingency

Minqiao Zheng, Zejun Yang*

(School of Mechatronic Engineering and Automation, Shanghai University, Shanghai 200444, China)

*Corresponding author at: Shanghai University, Shanghai 200444, China.

E-mail address: zyang@shu.edu.cn (Zejun Yang)

**Abstract:** Electrical power and computing power are highly interdependent especially due to the AI revolution. This dependency is largely underestimated during normal and recovery operation. To address this issue, an Electricity-Computing Interdependency System (ECIS) is developed in this paper, which first quantifies the interdependencies between ECIS and other critical infrastructures. The data center (DC) is modeled as a facility that transforms electrical power into computing power, and a recovery optimization strategy for ECIS is then proposed considering electrical power constraints, computing power supply, and infrastructure interdependencies. The case study demonstrates that the proposed strategy dynamically adjusts the DC electrical power quota according to system states, improving system resilience by more than 20% during power shortage.

## 1. Introduction

Computing power is a cornerstone of modern urban infrastructure. Critical infrastructures such as hospitals, financial institutions, urban water systems and communication networks all rely on stable computing services. Computing power supply interruptions directly affect the computing services and other daily activities of the public [1]. When extreme disasters damage the electrical power system, the resulting outages disrupt data centers (DCs) that provide computing support to other infrastructures. Given the strong interdependencies among these systems, such disruptions propagate through the network and lead to widespread cascading failures [2], ultimately reducing the overall functional level of the city.

While extensive research has focused on enhancing post-disaster electrical power system resilience and achieving critical load restoration under limited power resources [3], assessing urban resilience solely from the electrical power perspective can no longer reflect the operating characteristics of modern cities. This is because an increasing number of urban infrastructures rely heavily on external computing power resources [4]. For example, the intelligent diagnosis systems of hospitals [5] and the data processing platforms of financial institutions [6] require continuous computing services from DCs for their normal operation. Notably, some facilities (e.g., financial institutions) have limited electrical power consumption but rely heavily on computing power. In contrast, although DCs do not directly undertake urban public service functions, they consume electrical power to produce the computing power that other infrastructures depend on.

As shown in Fig. 1, based on infrastructure data from selected Chinese cities [7-15], this paper quantifies the differences in electrical power and computing power consumption across typical infrastructures. The analysis reveals a significant discrepancy between electrical power and computing power demands, both of which are critical for urban operations. Therefore, research on urban resilience systems with high dependence on computing power warrants more attention.

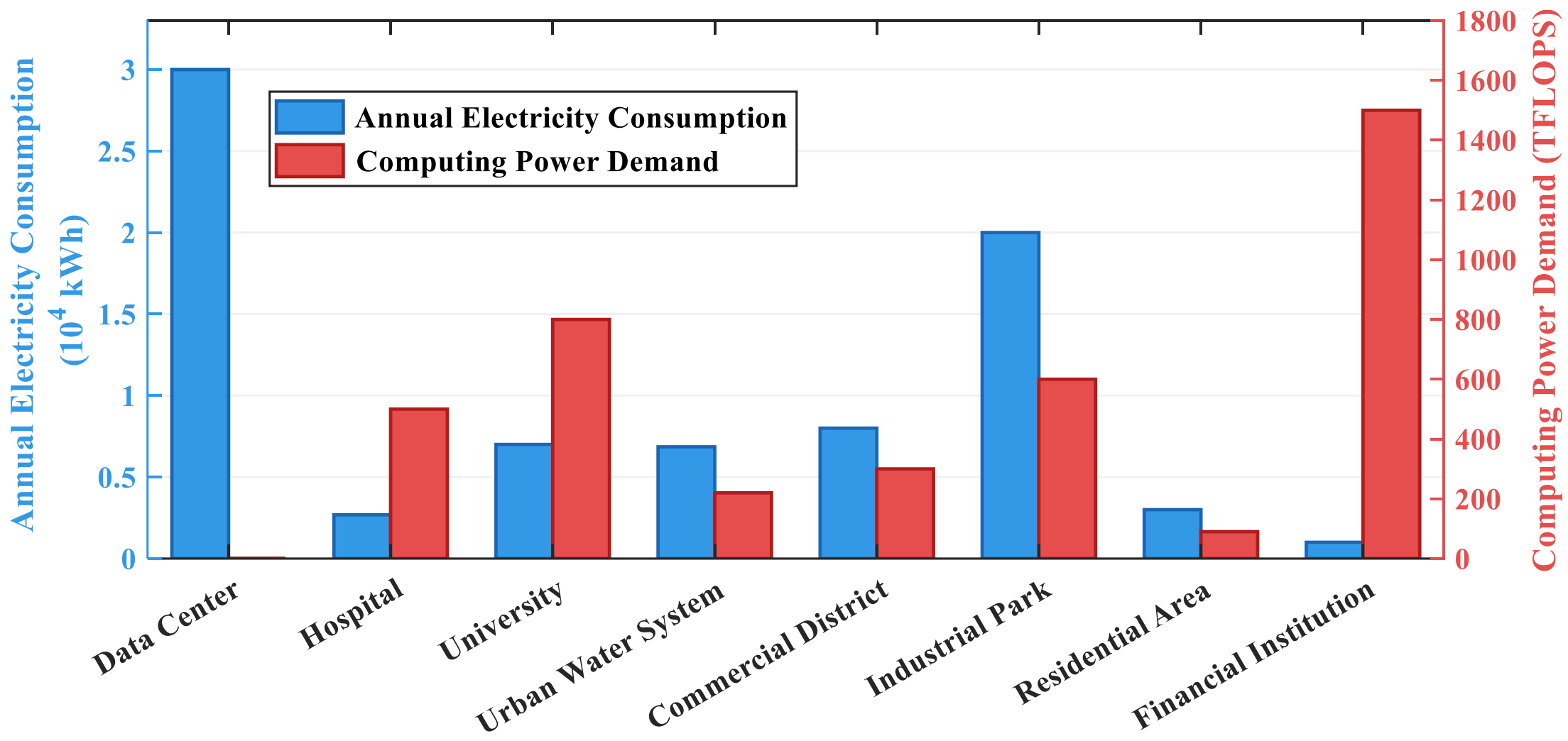

Fig. 1. Comparison of electrical power and computing power demand across different infrastructures.

Existing studies on electrical power system resilience have focused mainly on distribution network reconfiguration, mobile energy storage scheduling, repair crew optimization, power resource dispatch, and sequential restoration models [16-20]. Such studies typically prioritize electrical power supply to critical loads (e.g., hospitals and financial institutions) based on a predefined priority list to enhance system resilience [21, 22]. However, when a DC is also set as a load with a fixed priority level in post-disaster electrical power allocation, a core characteristic is overlooked: DCs are converting electrical power into computing power to support other infrastructures. For example, [23] optimizes the DC's own electrical power restoration through flexible substation regulation but still treats it as a conventional load without considering the role of computing supply for other urban infrastructures. Therefore, the traditional resilience model (if electrical power is on, then function is on) struggles to capture the computing power interdependencies in the AI revolution.

Although extensive research has been conducted on resource scheduling and collaborative optimization in computing power networks, covering aspects such as workload peak shaving for a single DC [24], geo-distributed task migration [25], and privacy-preserving multi-agent energy sharing [26], most of these studies implicitly assume sufficient electrical power resources and normal migration links. Their core objective remains limited to minimizing economic costs, without considering the constraints imposed by electrical power shortages under extreme disasters on computing power supply. In particular, existing literature lacks assessments of how insufficient computing power supply would cripple critical urban infrastructures and thereby reduce urban resilience [27]. Consequently, no unified framework has yet been established to describe the coordination among the computing power system, electrical power system, and urban infrastructures under extreme conditions.

To address the above issues, this paper proposes an Electricity-Computing Interdependency System (ECIS) framework, as shown in Fig. 2(a), to characterize the coupling relationships among computing power, electrical power, and other urban infrastructures under post-disaster scenarios. In the proposed framework, DC is modeled as a computing power production facility that consumes electrical power to provide computing services to different infrastructures. Based on this framework, this paper further develops a collaborative optimization model under post-disaster limited power conditions to enhance ECIS resilience (Fig. 2(b)).

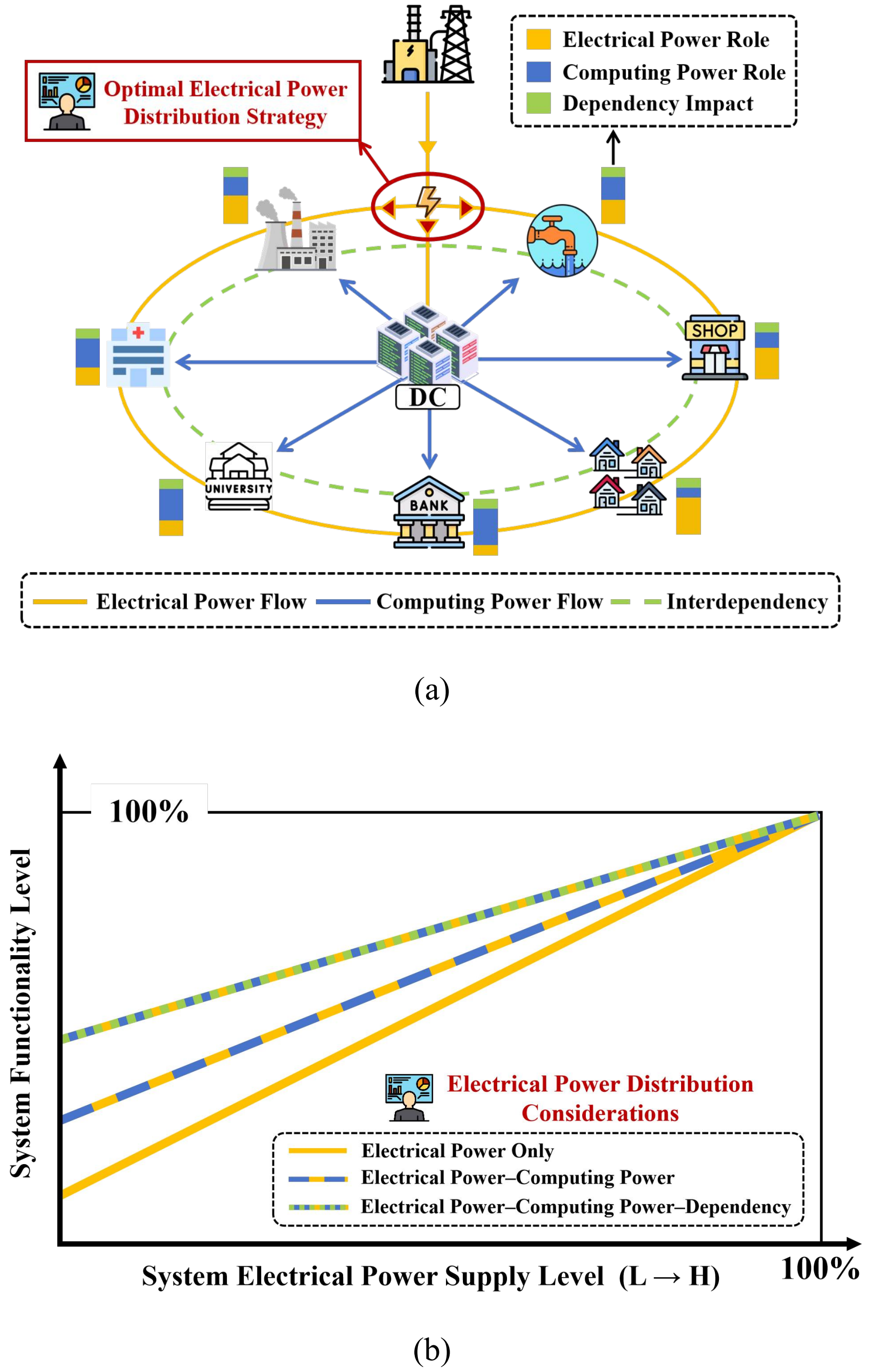

Fig. 2. (a) Framework of the ECIS; (b) post-disaster collaborative optimization performance.

The contributions of this paper are as follows:

1. An ECIS framework is proposed to characterize the coupling relationships among computing power, electrical power, and infrastructure interdependencies under post-disaster scenarios.
2. A collaborative resilience optimization strategy is proposed that considers computing power supply and infrastructure interdependencies, significantly improving system resilience compared to traditional predefined priority list restoration strategy.

# 2. System Description

## 2.1. Resilience Metric Definition

This section develops a functional model of urban systems that considers the coupling of "electrical power-computing power-infrastructure interdependencies". The ECIS resilience is described by the variation of system functional level with the electrical power supply ratio $k$ , which represents the ratio of available electrical power to normal at the current stage. Let the set of infrastructures be:

$$\mathcal{N} = \{1, 2, \ldots, N\} \tag{1}$$

where $N$ is the number of infrastructures. The data center is denoted as $dc \in \mathcal{N}$. It serves as both an electrical power load and a computing power supplier in the city, while the remaining facilities are computing power demand consumers. The overall system functional level is defined as:

$$F_{\text{sys}}(k) = \sum_{i \in \mathcal{N}} w_i F_i(k) \tag{2}$$

Here, $F_{\text{sys}}(k) \in [0,1]$ represents the functional level of the system at electrical power supply ratio $k$ (1 indicates full recovery, 0 indicates complete failure), $F_i(k)$ is the functional level of facility $i$ at that electrical power supply ratio, and $w_i \in [0,1]$ is the criticality weight of facility $i$, which measures its importance to overall urban operation.

$$\max R = \frac{1}{K} \sum_{k=1}^{K} F_{\text{sys}}(k) \tag{3}$$

This model uses discrete electrical power supply ratios $k$ to simulate the post-disaster recovery process and defines urban resilience as (3), where $K$ represents the number of recovery stages. Based on the system functional level defined above, the restoration optimization objective is to maximize the overall system resilience $R$ by optimizing the electrical power allocation at each recovery stage.

## 2.2. Facility Basic Function Model

The operation of modern urban infrastructures depends not only on electrical power resources but also on computing power resources provided by DCs. We define this operational function as the basic functional level $B_i$ and model it accordingly. Under post-disaster conditions, the actual electrical power supply received by an infrastructure is typically lower than the normal demand level. Therefore, the electrical power satisfaction level of facility $i$ is defined as:

$$P_i = \frac{P_i^{\text{alloc}}}{P_i^{\text{demand}}} \tag{4}$$

$$\sum_{i \in \mathcal{N}} P_i^{\text{alloc}} \leq P_{total}^{\text{avail}} \tag{5}$$

Here, $P_i^{\text{alloc}}$ is the actual allocated electrical power, $P_i^{\text{demand}}$ is the rated electrical power demand, and the sum of allocated electrical power for all infrastructures does not exceed the post-disaster available electrical power $P_{total}^{\text{avail}}$. Considering that actual infrastructures typically have a minimum operating electrical power requirement (e.g., operation of basic medical equipment in hospitals, pump station operation in urban water systems), a minimum operating electrical power threshold $P_i^{\text{min}}$ is defined. When $P_i < P_i^{\text{min}}$, the facility cannot produce effective functional output. For facilities that meet the minimum operating condition, the effective electrical power satisfaction level is defined as (6), where $\widetilde{P}_i \in [0,1]$ represents the effective power supply level after deducting the minimum operating requirement.

$$P_i = \frac{P_i - P_i^{\text{min}}}{1 - P_i^{\text{min}}} \tag{6}$$

This model focuses on the supporting role of DCs (as computing power supply facilities) in the recovery of external facility functions during post-disaster urban infrastructure restoration, rather than their internal operational efficiency or equipment reliability. Therefore, factors such as uninterruptible power supplies (UPS), diesel generators [28], cooling system energy consumption [29], and power usage effectiveness (PUE) [30] are not included in the model. Based on the above considerations, this paper adopts an "electrical power" to "computing power" mapping function to describe the overall relationship between the electrical power supply state of a DC and its computing power output capability:

$$C_{dc} = g\left(P_{dc}\right) \tag{7}$$

The electrical power satisfaction level of the DC is $P_{dc}$, and $g(\cdot)$ represents the nonlinear mapping relationship from electrical power to computing power. This model adopts a piecewise function to describe the computing power recovery process of the DC:

$$g(P_{dc}) = \begin{cases} 0, & P_{dc} < T_1 \\ g_1(P_{dc}), T_1 \le P_{dc} < T_2 \\ g_2(P_{dc}), T_2 \le P_{dc} < T_3 \\ g_3(P_{dc}), P_{dc} \ge T_3 \end{cases} \tag{8}$$

where $T_1$ is the startup threshold, $T_2$ is the ramping stage threshold, and $T_3$ is the saturation stage threshold. This model captures the computing power recovery process of a DC after a disaster: no effective computing power output (low electrical power availability), rapid increase in computing power (medium electrical power availability), and gradual saturation (high electrical power availability).

$$\sum_{i \in \mathcal{N} \setminus \mathrm{dc}} C_i^{\mathrm{alloc}} \le C_{dc} \tag{9}$$

$$C_i = \frac{C_i^{\mathrm{alloc}}}{C_i^{\mathrm{demand}}} \tag{10}$$

The computing power resources generated by the DC need to be allocated among different infrastructures. The computing power resource constraint is given by (9), where $C_i^{\mathrm{alloc}}$ is the computing power allocated to infrastructure $i$. Equation (10) defines the computing power satisfaction level of an infrastructure, and $C_i^{\mathrm{demand}}$ represents the computing power demand of infrastructure $i$.

The operational capability of an infrastructure is affected by both its electrical power supply state and computing power resources. Therefore, the basic functional level of facility $i$ is defined as:

$$B_i = \alpha_i P_i + \beta_i C_i \tag{11}$$

where $\alpha_i$ and $\beta_i$ are the weights of electrical power and computing power on facility functionality, respectively, satisfying $\alpha_i + \beta_i = 1$. A larger $\alpha_i$ indicates a higher dependence of the facility on electrical power, while a larger $\beta_i$ indicates a higher dependence on computing power. For example, residential areas mainly rely on electrical loads such as lighting, air conditioning, and elevators, thus having a higher electrical power weight; whereas financial institutions rely on real-time databases, cloud-based transactions, and digital financial services, thus having a higher computing power weight.

### 2.3. Infrastructure Dependency Propagation Model

When evaluating interdependencies among urban infrastructures, a typical "bucket effect" exists, where the overall system function is limited by its weakest link or most critical dependency. This concept is inspired by our previous work [2], which uses the Human Readable Table model to reveal weakest-link constraints among different emergency resources in power system resilience resource allocation. Building on this, the present work further extends it to urban infrastructure dependency scenarios.

$$A = \left[ a_{ij} \right] \tag{12}$$

$$D_i = \sum_{j=1}^{N} a_{ij} F_j \tag{13}$$

Equation (12) defines the dependency weight matrix, where $a_{ij} \geq 0$ represents the dependency intensity of facility $i$ on facility $j$. Accordingly, the dependency satisfaction level $D_i$ of facility $i$ is given by (13), where $F_j$ represents the actual functional level of facility $j$. This equation indicates that the dependency satisfaction level of a facility is jointly determined by its dependency weights on other facilities and the functional levels of those facilities. Therefore, the actual function of facility $i$ after considering dependencies is defined as:

$$F_i = B_i D_i \tag{14}$$

Here, $B_i$ represents its own basic functional level and $D_i$ represents its external dependency satisfaction level. This equation indicates that even if a facility's own electrical power supply and computing power are sufficient (i.e., $B_i$ is nominal), its actual function $F_i$ will still decrease if its critical dependent facilities fail. Since facility $i$ depends on the functions $F_j$ of other facilities $j$, and facility $j$ is in turn affected by the functions of other facilities, a direct analytical solution is not feasible. This model adopts an iterative approach to update facility functions step by step:

$$F_i^{(n+1)} = B_i \sum_j a_{ij} F_j^{(n)} \tag{15}$$

Here, $F_i^{(n)}$ represents the functional level of facility $i$ after the $n$-th iteration. In each iteration, the dependency satisfaction level of each facility $D_i^{(n)} = \sum_j a_{ij} F_j^{(n)}$, is computed from the functional values of all facilities in the previous iteration, and then multiplied by its own basic functional level $B_i$ to obtain the updated functional value for the current iteration. The convergence criterion for this iterative process is:

$$\max \left| F_i^{(n+1)} - F_i^{(n)} \right| < \varepsilon \tag{16}$$

That is, when the maximum change in functional values of all facilities between two consecutive iterations is less than the preset tolerance $\varepsilon$, the system is considered to have reached a steady state. For facilities with no external dependencies, their dependency satisfaction level is always 1(i.e., $F_i = B_i$).

# 3. Case Study

## 3.1. Parameter Settings

Based on the urban infrastructure electrical power-computing power demand discrepancy proposed in Section 1, this section validates the proposed ECIS resilience enhancement strategy. Extreme events are assumed to reduce the available electrical power below normal levels, and multiple scenarios are simulated with electrical power supply ratios ranging from 20% to 100% of the normal level. The electrical power and computing power demands of each facility are shown in Fig. 1. The DC converts electrical power into computing power and allocates it according to facility criticality levels. The functional level of each facility is jointly determined by electrical power and computing power, and is propagated iteratively through dependencies. To quantify functional importance, we classify infrastructures into three criticality levels: Level 1 ($w = 0.4$) for life-sustaining and economic cores (e.g., hospitals, financial institutions); Level 2 ($w = 0.35$) for public services (e.g., universities, urban water systems); and Level 3 ($w = 0.25$) includes commercial districts, industrial parks, and residential areas.

To validate the impact of interdependencies among the electrical power system, data center,

and other urban facilities on system resilience, Fig. 3 defines the interdependencies among different facilities and their dependence on electrical power. The thickness and color intensity of the arrows represent the dependency strength: thicker and darker arrows indicate higher dependency. Taking a financial institution as an example, its operations heavily depend on the DC, which is represented by a thick, dark arrow pointing to the DC. In contrast, its direct demand for electrical power in daily operations (e.g., lighting, office equipment, and air conditioning systems) is relatively small, so the arrow pointing to the power system is thinner and lighter.

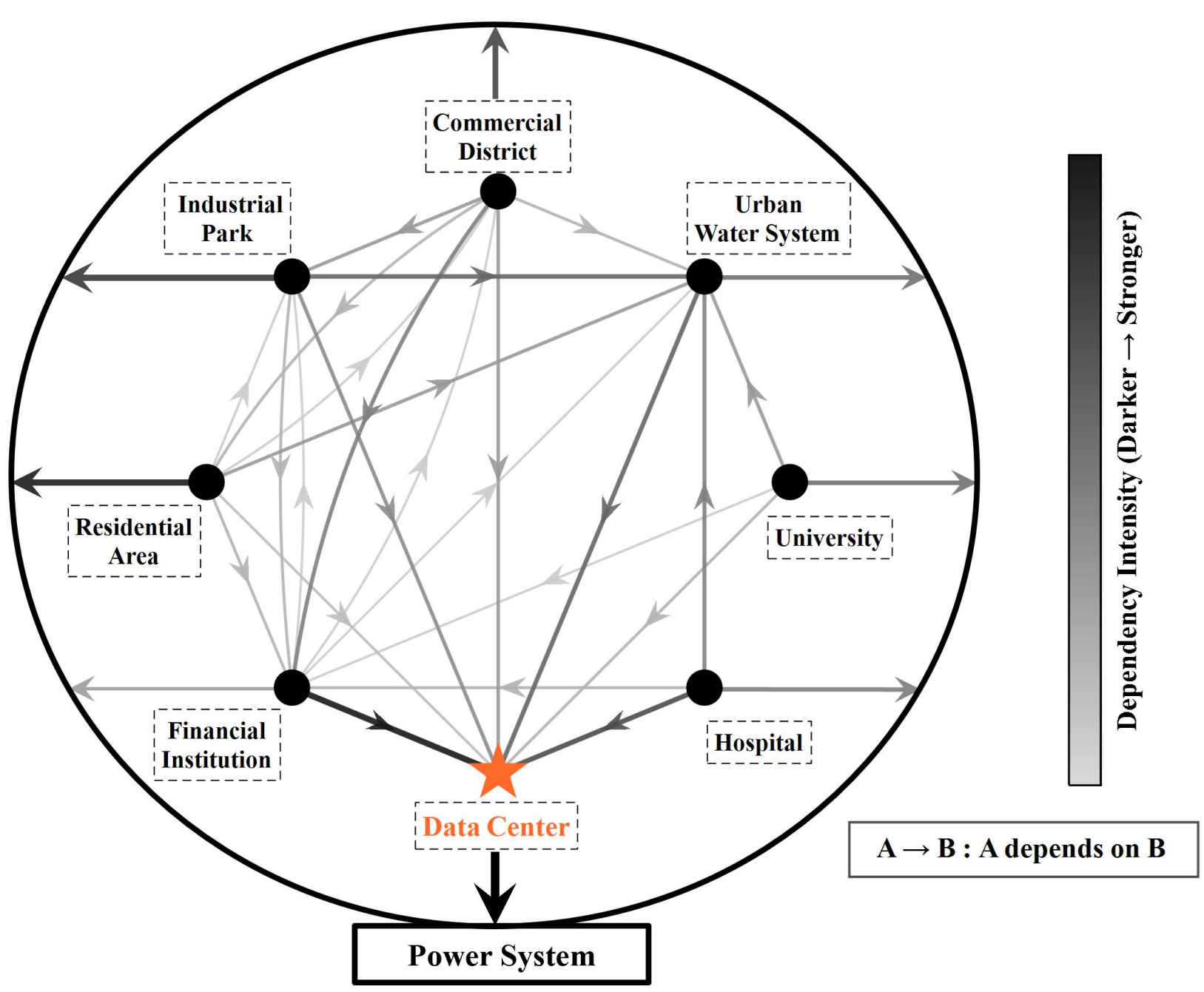

Fig. 3. Directed graph of interdependencies among urban infrastructures and their dependence on electrical power.

The optimization model constructed in this paper achieves joint decision-making of electrical power allocation, computing power supply, and infrastructure interdependencies under post-disaster limited electrical power conditions, which is essentially a nonlinear constrained optimization problem. Since neither the objective function nor the constraints are explicitly convex, this paper adopts a Genetic Algorithm for solution. Genetic Algorithm uses a population size of 100, a maximum of 300 generations, a crossover probability of 0.8, and an elite retention ratio of 5%.

As a computing power production facility, the DC adopts different power allocation schemes under different strategies. Two optimization strategies are tested to examine the impacts of different elements within the ECIS framework on system resilience:

- **Strategy A (predefined priority list restoration strategy):** only the contribution of electrical power resources to system resilience is considered. The DC is treated as a conventional load with a fixed priority level in the predefined list, and electrical power allocation is optimized solely according to a predefined priority list after a disaster.
- **Strategy B: (electrical-computing power interdependencies):** the computing power is considered as a resilience-enhancing factor and incorporated into electrical power allocation optimization. Furthermore, interdependencies among ECIS and other facilities are introduced to achieve muti-party collaborative optimization.

### 3.2. Strategy A: Predefined Priority List Restoration Strategy

Under limited post-disaster electrical power supply scenario, the electrical power is allocated solely based on a predefined priority list for each facility. The DC is modeled as a conventional load, without considering the impact of its computing power output on the functional levels of other facilities. Through comparing two scenarios in Strategy A where the DC is assigned a fixed criticality weight of 0.5 (high-weight case) versus 0.25 (low-weight case), the system functional level $F_{sys}$ (Fig. 4) and the composition of each facility's functional level $F_i$ (Figs. 5 and 6) are analyzed.

Notably, when the electrical power supply ratio is between 20% and 40%, the blue curve (corresponding to $w_{dc} = 0.5$, the high-priority case) shows a high value of $P_{dc}$ (the DC electrical power satisfaction level), yet the system functional level $F_{sys}$ remains at zero. When the electrical power supply ratio increases to 50%-80%, the system functional level ($F_{sys}$) of the blue curve exceeds that of the red curve. As electrical power recovery reaches 80%-100%, the DC also achieves a high electrical power satisfaction level ($P_{dc}$) under the low-criticality ($w_{dc} = 0.25$) scenario (Fig. 4(b)), and the system functional levels ($F_{sys}$) of the two scenarios converge.

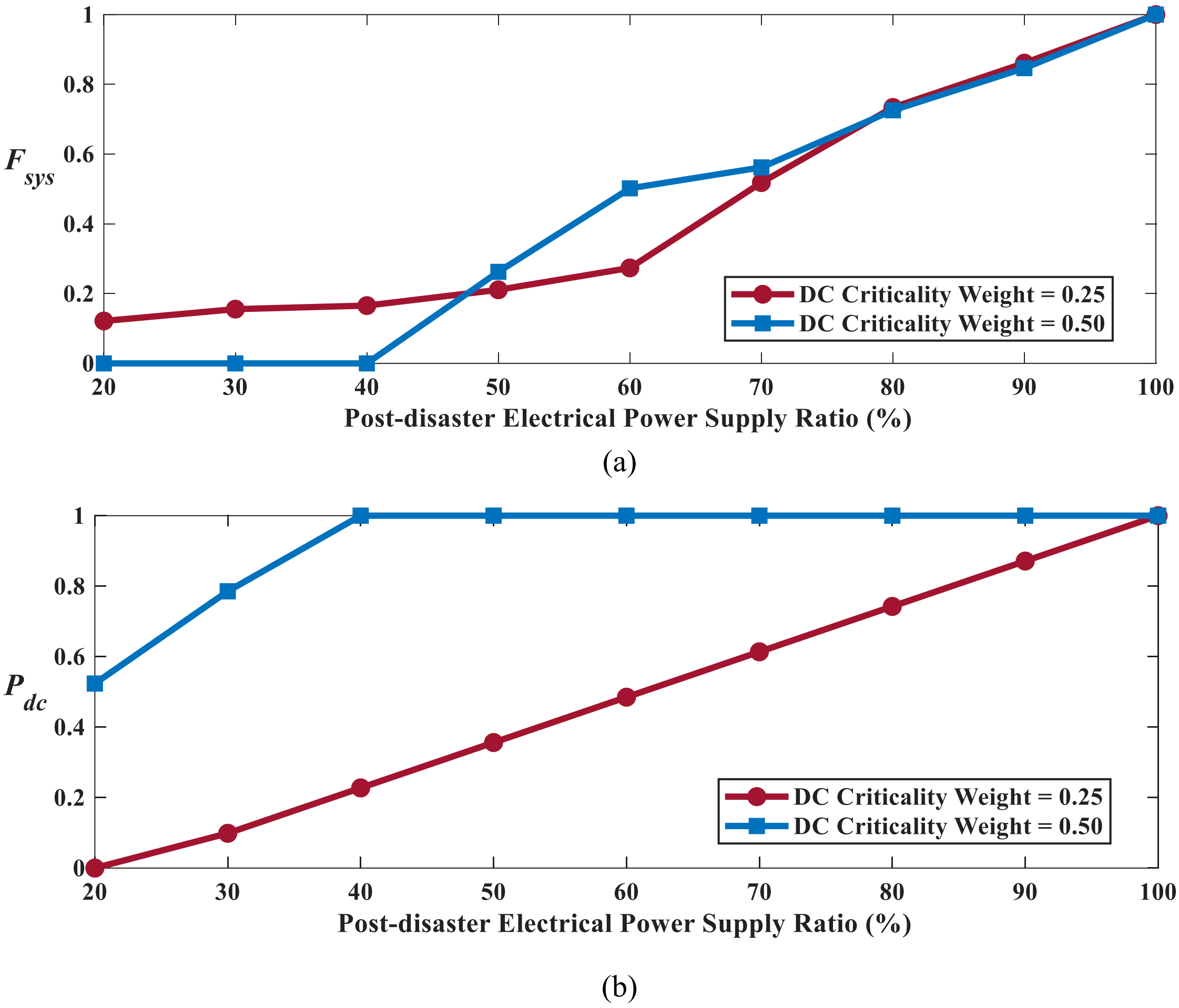

(a)

(b)

Fig. 4. Comparison between $w_{dc} = 0.5$ and $w_{dc} = 0.25$ under Strategy A: (a) system functional level $F_{sys}$; (b) DC electrical power satisfaction level $P_{dc}$.

Fig. 5 and 6 further explain the results in Fig. 4 by taking system electrical power supply ratios of 30% and 60% as examples, respectively. When the system electrical power supply ratio is 30% of the normal level (Fig. 5), if the DC is set as the most critical load ($w_{dc} = 0.5$), the system allocates all limited electrical power to the DC (upper subplot). In this case, the DC converts sufficient electrical power into computing power and maintains a high functional level. However, the remaining facilities receive no electrical power supply and cannot perform urban functions, resulting in a system functional level of zero (Fig. 4). Conversely, when the DC is set as a low-criticality load ($w_{dc} = 0.25$), the system prioritizes electrical power supply to other critical facilities (lower subplot). Although the DC receives no electrical power and thus cannot provide computing power, the other facilities can still achieve partial functional output under electrical power support.

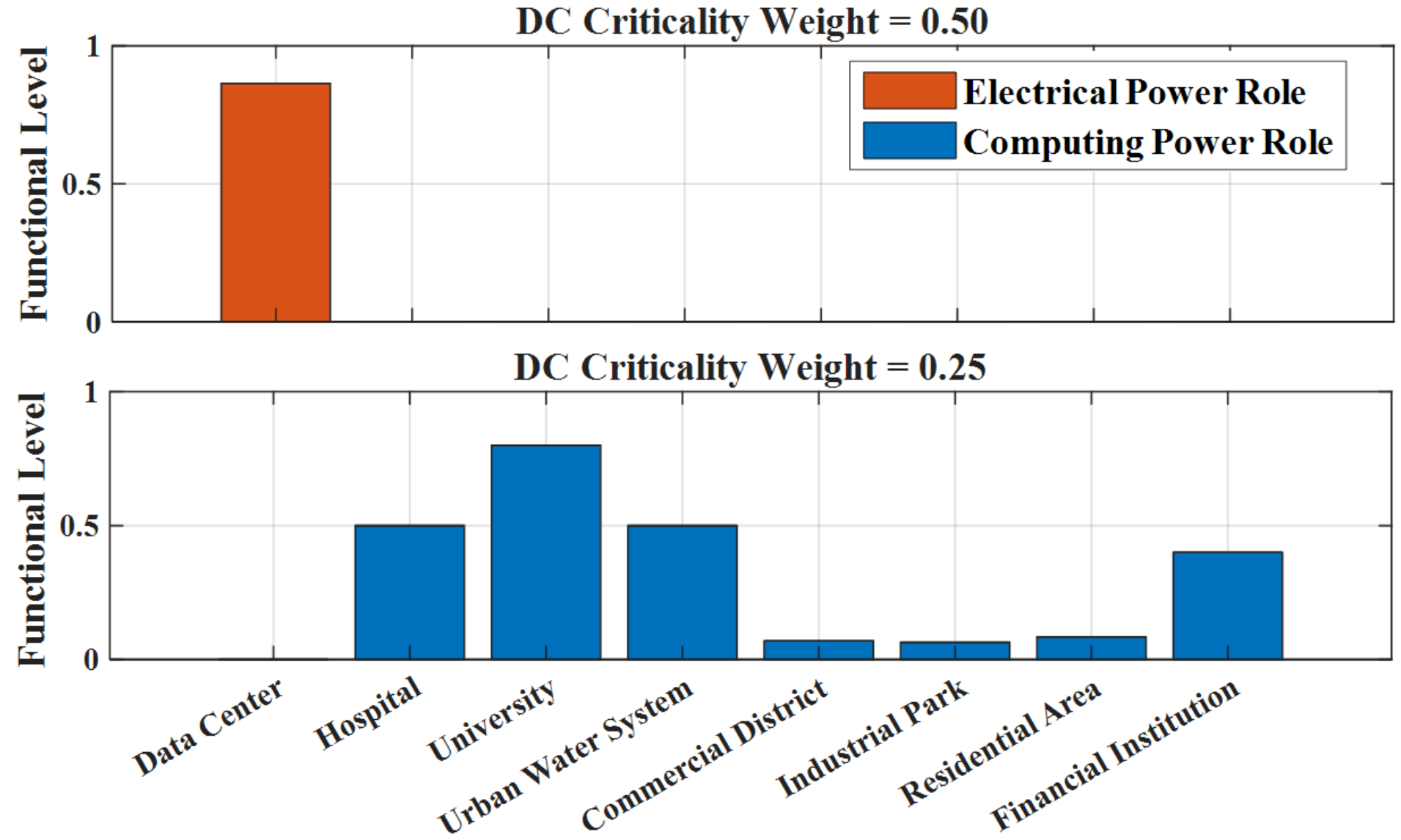

Fig. 5. Facility functional levels under $w_{dc} = 0.5$ and $w_{dc} = 0.25$ (system electrical power supply at 30% of the normal level).

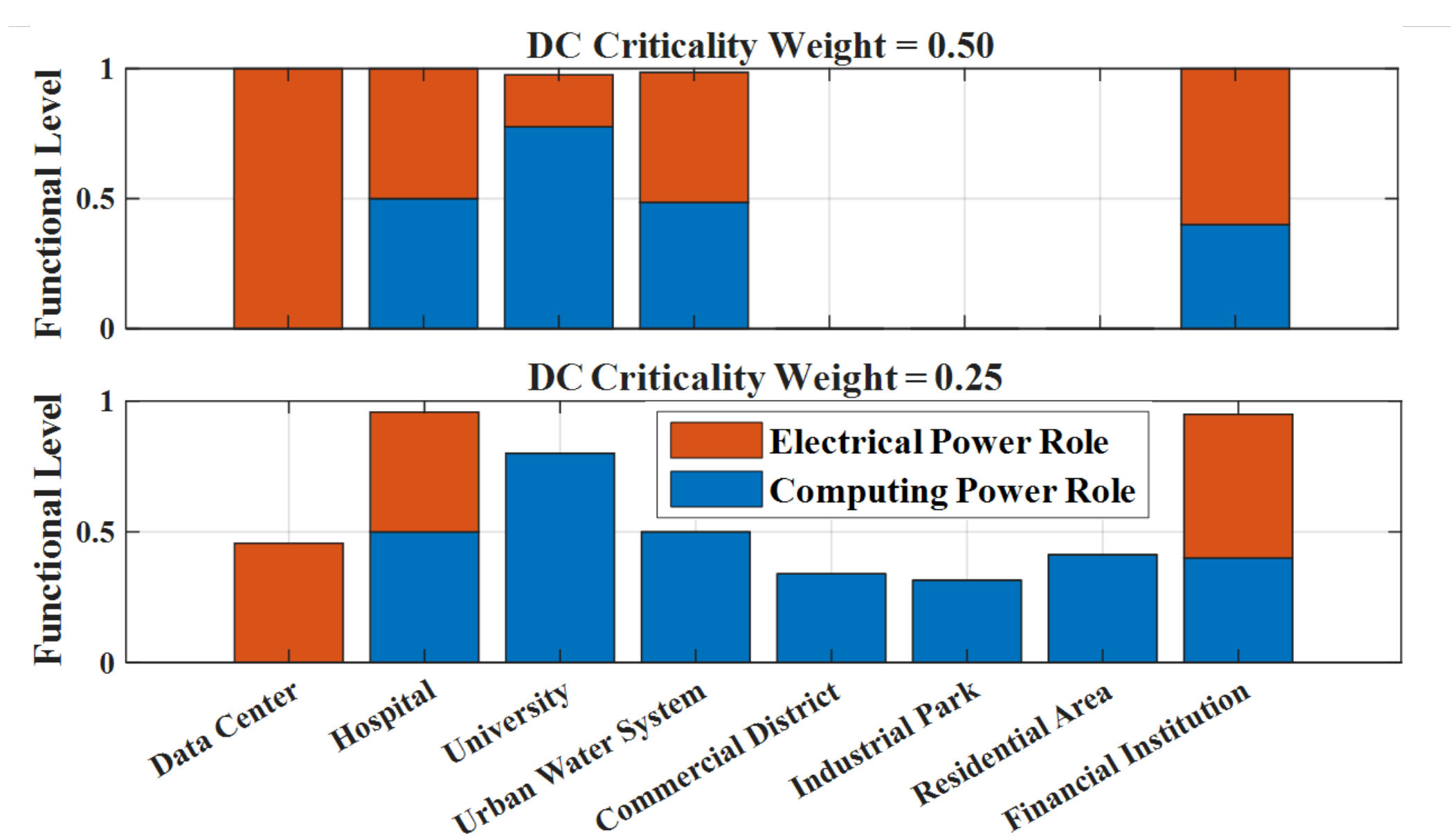

Fig. 6. Facility functional levels under DC weight = 0.5 and DC weight = 0.25 (system electrical power supply at 60% of the normal level).

When the electrical power supply ratio reaches 60% (Fig. 6), the results changes significantly. With the DC as the most critical load ($w_{dc} = 0.5$), Level 1, Level 2 facilities, and the DC all receive electrical power (upper subplot). The simultaneous provision of electrical and computing power allows these facilities to maintain a high functional level. With the DC as a low-criticality load ($w_{dc} = 0.25$), its electrical power priority drops, and power is distributed across all three levels of facilities (lower subplot). The DC then generates limited computing power, supporting only hospitals and financial institutions, while other facilities that rely on computing power for their functions cannot meet their needs. Thus, Level 2 facilities (university and urban water system) cannot achieve high functional levels. Since the system functional level calculation assigns greater weight to more critical facilities, setting the DC as the most critical load ($w_{dc} = 0.5$) at 60% electrical power supply (Fig. 4(a)) yields higher system functionality than setting it as low-criticality ($w_{dc} = 0.25$).

### 3.3. Strategy B (Electrical-Computing Power Interdependencies)

#### 3.3.1. Analysis of System Functional Levels

The results of Strategy A show that the performance implications of assigning the DC a fixed criticality weight of 0.5 versus 0.25 vary with the system electrical power supply ratio. Under severe power shortage, prioritizing electrical power supply to the DC leads to system functional collapse, whereas as the ratio increases, prioritizing the DC improves system functionality. This difference stems from the DC's electrical power-to-computing power conversion capability, which Strategy A overlooks. Therefore, Strategy B considers both the direct benefit of supplying electrical power to facilities and the indirect benefit of powering the DC to generate computing power, dynamically adjusting the DC's electrical power satisfaction level after a disaster.

As shown in Fig. 7(a), the green curve (Strategy B) performs best across the entire electrical power supply ratio range. Specifically: 1) In the low supply stage (20%-40%), Strategy B avoids the collapse of other facilities caused by over-prioritizing the DC (blue curve), maintaining a system functional level of approximately 0.2. 2) In the medium supply stage (50%-80%), it avoids the computing power shortage caused by neglecting the DC (red curve). 3) Even when electrical power supply recovers to above 80%, Strategy B still outperforms the two benchmark curves, demonstrating its advantage even when the electrical power shortage is minimal.

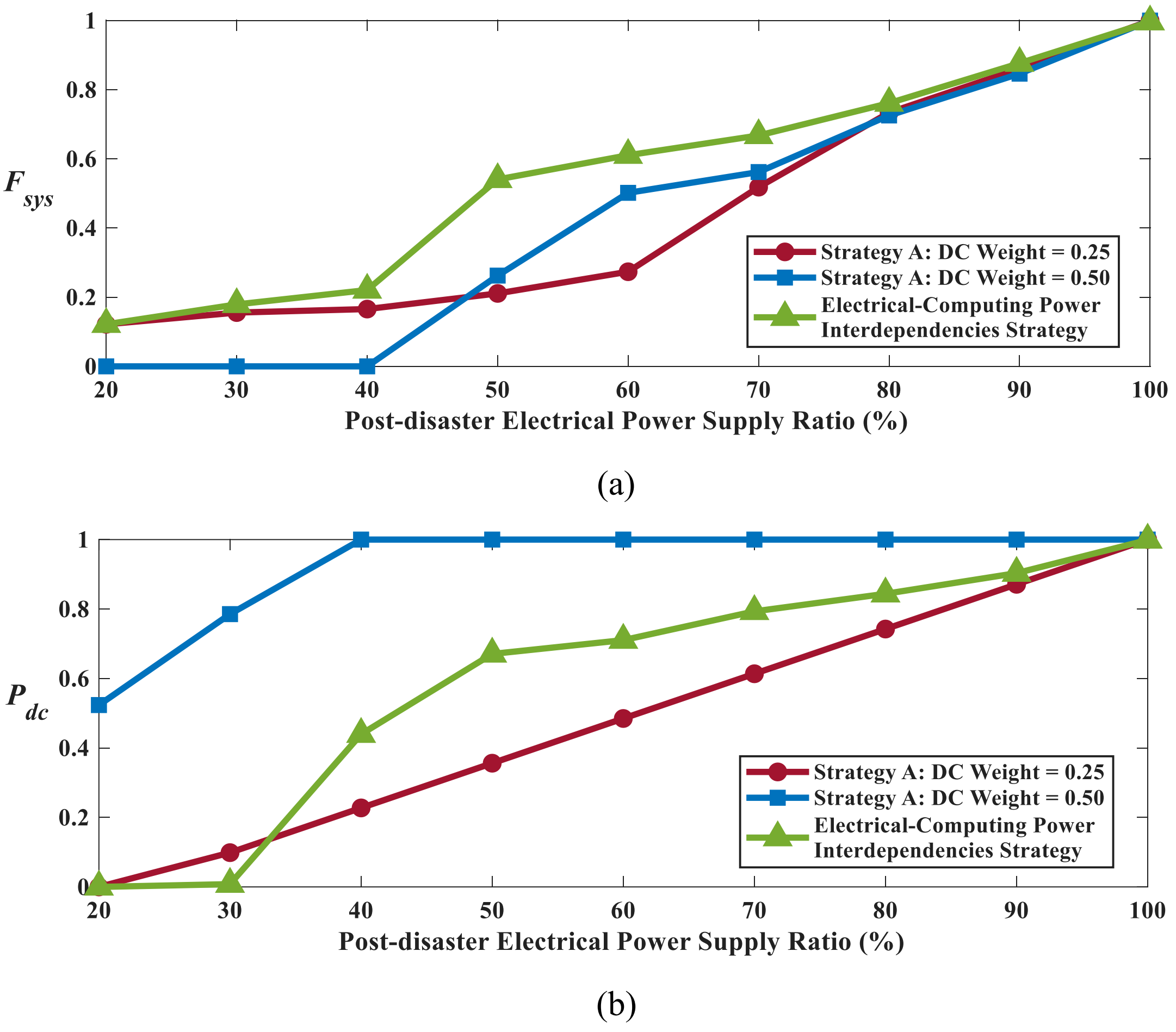

Fig. 7. Under the electrical-computing power interdependencies strategy: (a) system functional level $F_{sys}$; (b) DC electrical power satisfaction level $P_{dc}$.

Fig. 7(b) shows that Strategy B flexibly adjusts the DC's electrical power satisfaction level by comprehensively considering the available electrical power and system state. When the electrical power supply ratio is 20%-30%, the green curve value is zero, indicating that Strategy B temporarily sacrifices the DC to sustain fundamental urban functions. As electrical power supply gradually recovers, the electrical power demands of infrastructures are partially met, and the system begins to allocate more electrical power to the DC. At this point, the benefit of computing power to other facilities becomes evident, prompting the strategy to dynamically adjust and increase electrical power support to the DC.

### 3.3.2. Impact Mechanism of Interdependencies on Strategy B

To demonstrate the importance of considering interdependencies among infrastructures in electrical power allocation optimization, this session compares system resilience under two scenarios ("with dependencies" vs. "without dependencies") to reveal the impact of interdependencies on post-disaster recovery.

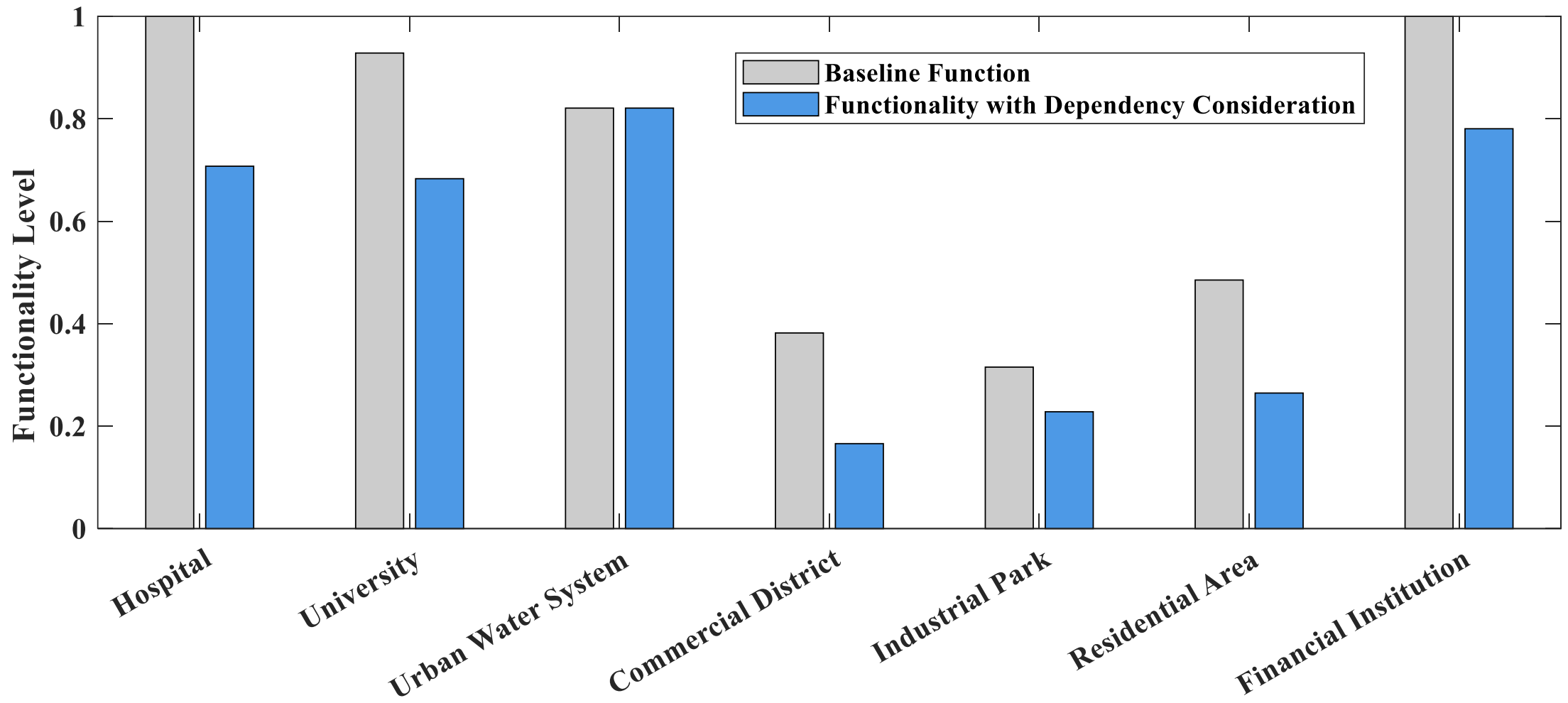

Fig. 8. Comparison of functional levels with and without considering infrastructure interdependencies.

Fig. 8 illustrates the impact of infrastructure interdependencies on post-disaster functional levels. The baseline function (gray bars) considers only the direct effects of electrical power and computing power, ignoring the impact of interdependencies. Under this setting, for highly critical facilities (hospitals and financial institutions), their electrical power and computing power needs are prioritized, allowing them to reach 100% if dependencies on other facilities are ignored. However, introducing the dependencies shown in Fig. 3 reduces their actual functionality due to dependency constraints. Notably, as specified in Fig. 3, the urban water system depends only on the electrical power system and DC. Since the baseline function already accounts for both electrical power and computing power factors, its functional level remains unchanged under the "with dependencies consideration" scenario. Strategy B incorporates infrastructure interdependencies as a key constraint in electrical power allocation optimization. Table I quantitatively compares the system resilience values of Strategy B with and without considering dependencies. The results show that considering dependencies improves system resilience by 9.5% (from 0.505 to 0.553) compared to the case without dependencies (Strategy B without Considering Dependencies).

Table I. Performance comparison between electrical-computing power interdependencies strategy and predefined priority list restoration strategy.

| **Post-disaster Power Supply Ratio** | **40%** | **60%** | **80%** | **Resilience** |
|---|---|---|---|---|
| Strategy A ($w_{dc} = 0.25$) | 0.166 | 0.273 | 0.734 | **0.449** |
| Strategy A ($w_{dc} = 0.50$) | 0 | 0.502 | 0.726 | **0.433** |
| Strategy B without Considering Dependencies | 0.180 | 0.572 | 0.725 | **0.505** |
| Strategy B: Electrical-Computing Power Interdependencies | 0.220 | 0.610 | 0.761 | **0.553** |

Table I also compares the performance differences between Strategy B and the two predefined priority list restoration strategy (Strategy A). The results show that Strategy B achieves a resilience value of 0.553, which is 23.2% and 27.7% higher than those of Strategy A with $w_{dc} = 0.25$ (0.449) and with $w_{dc} = 0.5$ (0.433), respectively. These results confirm that collaborative optimization incorporating computing power and interdependencies overcomes the limitations of the predefined priority list restoration strategy and maximizes system resilience. Meanwhile, ignoring interdependencies among facilities leads the model to overestimate the actual functional levels of post-disaster urban infrastructures, resulting in an overly optimistic resilience assessment by the electrical power sector.

# 4. Conclusion

This paper proposed an Electricity-Computing Interdependency System framework for post-disaster restoration to characterize the coupling relationships among computing power, electrical power, and infrastructure interdependencies. The case study results demonstrated that the traditional restoration strategy based on a predefined priority list struggles to adapt to the widespread discrepancy between electrical power demand and computing power demand in modern cities. However, the proposed electrical-computing power interdependencies strategy can dynamically adjust the DC's electrical power quota according to system states, thereby improving system resilience by more than 20% under post-disaster power shortage.

This study revealed that in the context of the rapid development of artificial intelligence and digital infrastructures, computing power has become a critical resource affecting urban resilience. Assessing post-disaster recovery capability solely from the electrical power dimension can no longer accurately reflect urban operation mechanisms. Future work will further incorporate real distribution network topologies, sequential restoration processes, and multi-data-center coordination scenarios to investigate the resilience of ECIS in large-scale urban environments.